\documentclass[11pt]{amsart}

\usepackage{amssymb} 
\usepackage[mathscr]{eucal} 
\usepackage{xypic}   
\usepackage{amsmath} 
\usepackage{epsfig}
\usepackage{amscd}

\newtheorem{TEO}{Theorem}[section]
\newtheorem{PROP}[TEO]{Proposition}
\newtheorem{LEM}[TEO]{Lemma}

\newtheorem{COR}[TEO]{Corollary}
\newtheorem{REM}[TEO]{Remark}

\newcommand\dual{\mathrel{\raise3pt\hbox{$\underline{\mathrm{\thinspace d
\thinspace}}$}}}

\newcommand\proj{\mathbb P}
\newcommand\Z{\mathbb Z}

\newcommand\R{\mathbb R} 
\newcommand\Co{\mathbb C}

\newcommand\Ag{{A}_g^{(n)}}
\newcommand\Mg{{M}_g^{(n)}}

\def\Z{{\mathbb Z}}
\def\R{{\mathbb R}}
\def\Z{{\mathbb Z}}

\begin{document}


\title{Siegel metric and curvature of the moduli space of curves}
\author{Elisabetta Colombo\\ Paola Frediani}
\date{}
%
%

\footnote{
The present research took place in the framework of the
PRIN 2005 MIUR: "Spazi dei moduli e teoria di Lie" and PRIN 2006 of MIUR "Geometry on algebraic varieties". 

 }

\maketitle

\begin{abstract}
We study the curvature of the moduli space  ${M_g}$ of curves of
genus $g$ with the Siegel metric induced by the period map $j:{
M_g}\rightarrow {A_g}$. We give an explicit formula for the
holomorphic sectional curvature of ${M_g}$ along
 a Schiffer variation $\xi_P$, for $P$ a point  on the curve
 $X$, in terms of the holomorphic sectional curvature of ${A_g}$ and the second Gaussian map.
 Finally we extend the
K\"ahler form of the Siegel metric as a closed current on
$\overline{M}_g$ and we determine its cohomology class as a
multiple of $\lambda$.
\end{abstract}

\section{Introduction}\label{sec:1}

 During the last thirty years, some
natural  metrics on the moduli space of genus $g$ curves $M_g$
have been extensively studied. Many of these metrics come from
metrics on the Teichm\"uller space of which the moduli space of
curves is the quotient by the mapping class group. One of these is
the Weil-Petersson metric $\omega_{WP}$. It was introduced by Weil
and is known to be K\"ahler, to have non-positive curvature
operator and negative Ricci curvature, and to be geodesically
convex. S. A. Wolpert showed that both its holomorphic sectional
curvature and Ricci curvature have negative (genus dependent)
upper bounds (but no lower bounds do exist). Moreover it is not
complete (\cite{wol1}). The other canonical metrics, namely the
Teichm\"uller metric (or the Kobayashi metric), the Caratheodory
metric, the K\"ahler-Einstein metric, the induced Bergman metric,
the McMullen metric,  are complete. Recently Liu, Sun and Yau
(\cite{lsy1}, \cite{lsy2}) showed their equivalence on $M_g$ and
the equivalence with the Ricci metric and the perturbed Ricci
metric introduced by them.
 The K\"ahler form of the WP metric has been extended by Masur (\cite{ma})  as a closed current on the
 Deligne-Mumford compactification $\overline{M}_g$ of $M_g$.
 Wolpert (\cite{wol}) determined its cohomology class in terms of
 the first Chern class of the Hodge bundle $\lambda$ and the classes of the boundary.

Let $A_g$ be the moduli space of principally polarized abelian
varieties of dimension $g$ and let $j: M_g\rightarrow A_g$ be the
period map sending a curve to its jacobian. It is an interesting
and classical problem to understand the geometry of the image of
$M_g$ in $A_g$.

On $A_g$ there is a natural metric coming from the unique
$Sp(2g,\R)$ invariant metric on the Siegel space $H_g\simeq
Sp(2g,\R)/U(g)$ of which $A_g$ is the quotient by $Sp(2g,\Z)$.
 The purpose of this paper is to study the metric on $M_g$ induced by this metric through the period
map, which we call the Siegel metric. In \cite{cpt} an explicit
expression for the second fundamental form of the immersion $j$ is
given and it is proven that the second fundamental form lifts the
second Gaussian map $\mu_2:I_2(K_X)\rightarrow H^0(X,4K_X)$, as
stated in an unpublished paper of Green-Griffiths (cf.
\cite{green}).

Here we use it to compute the curvature of the Siegel metric. In
particular we give an explicit formula for the holomorphic
sectional curvature of $M_g$ along the
 a Schiffer variation $\xi_P$, for $P$ a point  on the curve
 $X$, in terms of the holomorphic sectional curvature of $
 A_g$ and the second Gaussian map $\mu_2:I_2(K_X)\rightarrow
H^0(X,4K_X)$.

Finally we give some properties of the holomorphic sectional
curvature of $M_g$, using  results of \cite{ cfW}. In particular
along a Schiffer variation $\xi_P$ the holomorphic sectional
curvature $H(\xi_P)$ of $ M_g$ is strictly smaller than the
holomorphic sectional curvature of $A_g$ unless $P$ is either a
Weierstrass point of a hyperelliptic curve or a ramification point
of the $g^1_3$ on a trigonal curve. In these last cases
$H(\xi_P)=-1.$

Furthermore we study the asymptotic behaviour of  the K\"ahler
form of the Siegel metric on $M_g$ showing that it extends as a
closed current to $\overline{M}_g$, hence it defines a cohomology
class in $H^2(\overline{M}_g, \Co)$ which we compute to be $\pi
\lambda$.

In all what we have stated, we have considered $M_g$ as an
orbifold. In the paper we make all computations using the covering
of $M_g$ given by the moduli space of curves with level $n\geq 3$
structures $\Mg$ and the moduli space $\Ag$ of principally
polarized abelian varieties with level $n$ structures.

In fact $\Mg$ and $\Ag$ are smooth and by Local Torelli theorem
proven in \cite{os} we know that the period map $j^{(n)}: \Mg
\rightarrow \Ag$ is a two to one immersion outside the
hyperelliptic locus and it is an injective immersion if we
restrict to the hyperelliptic locus.

The paper is organized as follows: in Section 2 we define the
Siegel metric and compute it on the tangent directions given by
the Schiffer variations (Lemma \ref{scalprod}).
 In Section 3 we give the expression of the curvature of the
Siegel metric on $ \Ag$ restricted to the Schiffer variations.
Then, we show that the second fundamental form of the immersion of
$\Mg$ in $\Ag$ is non zero at any non hyperelliptic curve and we
exhibit a formula for the curvature of $\Mg$ (Thm.\ref{formula}).
Finally we write the holomorphic sectional curvature of $\Mg$
along a Schiffer variation $\xi_P$, using the second Gaussian map.
 In Section 4 we give some applications of results of \cite{cfW} to the holomorphic sectional curvature of
$\Mg$. In Section 5 we study in particular the hyperelliptic locus
$HE_g$ and we show that the second fundamental form of $HE_g$ in
$\Ag$ is non zero at any point. In Section 6 we extend the
K\"ahler form of the Siegel metric as a closed current on
$\overline{M}_g$ and we determine its cohomology class
(\ref{class}).

{\bf Acknowledgments.} The authors thank Gilberto Bini and Pietro Pirola for several
fruitful suggestions and discussions on the subject.

\section{The Siegel metric}\label{sec:2}

We introduce some notations.  Let $M_g$, resp. $M_g^{(n)}$ be the
moduli space of smooth genus $g$ curves, resp. of smooth genus $g$
curves with a fixed $n$-level structure. Denote by $T_g$ the
Teichmuller space and by $\Gamma_g$ the mapping class group acting
on $T_g$ with quotient $M_g$. Let $K(n):= ker (\Gamma_g\rightarrow
Sp(2g,\Z/n\Z))$ and recall that $M_g^{(n)}$ is the quotient of
$T_g$ by the action of $K(n)$. Moreover, let $K:= ker
(\Gamma_g\rightarrow Sp(2g,\Z))$ be the Torelli group and define
the Torelli space $Tor_g$ as the quotient of $T_g$ by the action
of the Torelli group.

Let $A_g$, resp.\ $A_g^{(n)}$ be the moduli space of
$g$-dimensional principally polarized Abelian varieties, resp.\ of
$g$-dimensional principally polarized Abelian varieties with a
$n$-level structure. Denote by $H_g:=\{Z\in M(g,\Co)\ |\ Z=^tZ,
ImZ>0\}$ the Siegel space so that $A_g$ is the quotient of $H_g$
by the action of $Sp(2g,\Z)$ and $A_g^{(n)}$ is the quotient of
$H_g$ by $ker(Sp(2g,\Z)\rightarrow Sp(2g,\Z/n\Z))$. Denote by
$j^{Tor}$, $j$ and $j^{(n)}$ the period maps which send a curve to
its jacobian. We have the following diagram
$$
\begin{array}{ccc}
T_g & & \\
\downarrow & & \\
 Tor_g & \stackrel{j^{Tor}}\rightarrow& H_g \\
 \downarrow & & \downarrow\\
 M_g^{(n)}& \stackrel{j^{(n)}}\rightarrow& A_g^{(n)}  \\
  \downarrow & & \downarrow \\
 M_g
 & \stackrel{j}\rightarrow& A_g \\
\end{array}
$$
 The Torelli theorem states that $j$ is injective, while $j^{tor}$ and $j^{(n)}$ are two to one on the
 image and ramified over the hyperelliptic locus.
In fact multiplication by $-1$ in $H^1(X,\Z)=H^1(JX,\Z)$, where
$JX$ is the Jacobian of the curve $X$, is induced by an
automorphism of abelian
 varieties but not by an automorphism of non hyperelliptic curves.
 Local Torelli Theorem says that outside the hyperelliptic
 locus and restricted to the hyperelliptic locus the period map is an immersion (cf. \cite{os}).
 From now on we shall work on $\Mg$ and $\Ag$, with $n\geq 3$, since they are
 smooth,  everything works in the same way on  $M_g$ and $A_g$ but in the
 orbifold context.

 We will now define the Siegel metric.

The Siegel space $H_g$ is a homogeneous space and it can be seen
as the quotient $Sp(2g, \R)/U(g)$. We call the unique (up to
scalar) invariant metric the Siegel metric.

OnLet $F$ be the homogeneous vector bundle on $H_g$ associated to
the standard $g$-dimensional representation of $U(g, \Co)$. The
Hodge metric $h$ on $F$ is the only (up to multiplication by
scalars) invariant metric on the homogeneous bundle $F$. Moreover
through the identification
$$\Omega^1_{H_g}\simeq S^2F$$ the Hodge metric on $F$ defines the
Siegel metric on $H_g$ .

The Siegel metric on $H_g$ defines a metric on $\Ag$ and $A_g$
and, through the period map,  an induced metric on $\Mg$ and $M_g$
outside the hyperelliptic locus, and on the hyperelliptic locus
itself. We call all these metrics the Siegel metrics.

These metrics can be described in terms of polarized variation of
Hodge structures. More precisely, on $\Ag$ we have the universal
family $\phi : {\mathcal A} \rightarrow \Ag$, and the polarized
variation of Hodge structures associated to the local system
$R^1\phi_*\Z$. The associated Hodge bundle ${\mathcal F}^1$ can be
identified with $\phi_* (\Omega^1_{{\mathcal A}|\Ag})$, where
$\Omega^1_{{\mathcal A}|\Ag}$ is the sheaf of relative holomorphic
one forms. The polarization induces a Hermitian metric on
$R^1\phi_*\Co$ and on ${\mathcal F}^1$, which we call the Hodge
metric. In fact the pullback of ${\mathcal F}^1$ on $H_g$ is the
bundle $F$ and the pullback of the metric is the Hodge metric on
$F$. Hence the Siegel metric is induced by the Hodge metric
through the identification $S^2 {\mathcal F}^1 \cong
\Omega^1_{\Ag}$.

On $\Mg$ we have the universal family $\psi: {\mathcal C}
\rightarrow {M_g^{(n)}}$ with induced relative dualizing sheaf
$K_{{\mathcal C}|{M_g^{(n)}}}$. The local system $R^1\psi_*\Z$
coincides with the pullback of $R^1\phi_*\Z$ through the period
map: at a point $ [X] \in \Mg$, we have $H^1(X, \Z) \cong H^1(JX,
\Z)$.  The non-degenerate Hermitian product on $H^1(X,{\mathbb
C})$, defined by the polarization is the following: for any
$[\eta], [\xi] \in H^1(X,{\mathbb C})$, we have
$$ \langle [\eta], [\xi] \rangle = i \int_{X} \eta \wedge
\overline{\xi}.
$$
 The Hodge bundle can be identified with $\psi_*({K}_{{\mathcal C}|\Mg})$, and the corresponding
 Hodge metric yields a metric on $S^2 {\mathcal F}^1 \cong {j^{(n)}}^*\Omega^1_{\Ag}$,
 hence on $ {j^{(n)}}^*{\mathcal T}_{\Ag}$, and by restriction the Siegel metric on ${\mathcal T}_{\Mg}$.

We finally observe that for the sake of simplicity we defined the
Siegel metric on the fine moduli space $\Mg$, but we also have a
Siegel metric on $M_g$ viewed as an orbifold.

\subsection{An explicit formula}
We shall now give an explicit formula for the Siegel metric on
$\Mg$ at a point $[X] \in \Mg$ in terms of the basis of $H^1(T_X)$
given by Schiffer variations $\xi_P$, for a set of $3g-3$ general
points on $X$.

Now, we briefly recall the definition of $\xi_P$. Consider the exact
sequence
$$
0 \rightarrow T_X \rightarrow T_X(P) \rightarrow T_X(P)_{|P} \rightarrow 0.
$$

Notice that $H^0(T_X(P)_{|P}) \cong {\Co}$. If we denote the
coboundary map by $\delta: H^0(T_X(P)_{|P}) \rightarrow H^1(T_X)$,
we have $dim(Im(\delta)) = 1$. Any no zero element $\xi_P$ in
$Im(\delta)$ is called a Schiffer variation. Let us a choose a
local coordinate $z$ in a neighborhood of $P$. Under the Dolbeault
isomorphism $H^1(T_X) \cong H^{0,1}(T_X)$, it is represented by
the form
$$
\theta_P = \frac{1}{z} \overline{\partial}b_P \otimes
\frac{\partial}{\partial z},
$$ where $b_P$ is a bump function around $P$.
Notice that if we choose $b_P$ to be one in a neighborhood of $P$
for this choice of local coordinate $z$, $\xi_P$ depends only on
the choice of $z$.  In what follows, we need to express $\xi_P$ in
terms of a basis of $S^2(H^0(K_X)^*)$, through the inclusion of
$H^1(T_X)$ in $S^2(H^0(K_X)^*)$.

Fix an orthonormal basis $\{\omega_i\}_{i=1, \ldots, g}$ of
$H^0(K_X)$. Choose a local coordinate $z$ around a point $P \in X$
and write $\omega_j = f_j(z) \, dz$. Since $H^0(K_X) \cong
H^{1,0}(X)$, the set $\{\overline{\omega}_i\}_{i=1, \ldots, g}$ is
a basis of $H^{0,1}(X)$. This set can be viewed as the dual basis
of $\{\omega_i\}$, where the non-degenerate pairing is given by
$$
\overline{\omega}_i(\omega_j) = i \int_X \omega_i \wedge
\overline{\omega}_j =\langle \omega_i, \omega_j \rangle = \delta_{ij}.
$$

Observe that we have
\begin{equation}
\label{campro}
\langle \overline{\omega}_i, \overline{\omega}_j \rangle = i \int_X
\overline{\omega}_i \wedge \omega_j = - i \int_X \omega_j \wedge
\overline{\omega}_i = - \delta_{ij}.
\end{equation}

\begin{LEM}
\label{csip}  For a choice of a local coordinate $z$ at $P$, we
have
\begin{equation}
\label{stoka}
\langle\xi_P(\omega_i),\overline{\omega_j}\rangle =
-2 \pi f_i(P) f_j(P),
\end{equation}
 Hence in $S^2(H^0(K_X)^*))\cong S^2(H^{0,1}(X))$ it holds
$$
\xi_P = 2\pi \sum_{i,j}^g f_i(P)f_j(P) (\overline{\omega}_i \odot
\overline{\omega}_j).
$$
\end{LEM}

\proof Since $ \xi_P = \sum_{i =1}^g \xi_P(\omega_i) \odot
\overline{\omega}_i$ and
\begin{equation}
\label{xdopo}
\xi_P(\omega_i)= \sum_j -\langle \xi_P(\omega_i),
\overline{\omega}_j \rangle \overline{\omega}_j,
\end{equation}
we get
$$
\xi_P = \sum_{i,j}  -\langle \xi_P(\omega_i),
\overline{\omega}_j \rangle \left( \overline{\omega}_j \odot \overline{\omega}_i\right).
$$

By definition of $\xi_P$, the element $\xi_P(\omega_i) \in H^{0,1}(X)$
is represented by the $(0,1)$-form
$$ (\frac{1}{z}\overline{\partial}b_P \otimes \frac{\partial}{\partial
z}) \rfloor(f_i(z) dz) = \frac{1}{z} \overline{\partial} b_P f_i(z).
$$

Let $C$ be a small circle around $P$ such that $b_P \equiv 1$ on
$C$. Then the lemma follows by Stokes and Cauchy Theorems:
$$
\langle\xi_P(\omega_i),\overline{\omega_j}\rangle =
i \int_{X} \frac{\overline{\partial} b_P}{z-z(P)} f_i(z) \wedge f_j(z) dz =
i \int_X \overline{\partial} \left(\frac{b_Pf_i(z) f_j(z)}{z-z(P)} \right) \wedge dz =
$$
$$
i \int_X d \left(\frac{b_Pf_i(z) f_j(z)}{z-z(P)}dz \right) = i
\int_C \frac{f_i(z) f_j(z)}{z - z(P)} dz=-2\pi f_i(P)f_j(P)
$$

\qed

\begin{LEM}
\label{scalprod}
The scalar product of the two Schiffer variations $\xi_P$, $\xi_{P'}$ has the
following form:

$$\langle \xi_P, \xi_{P'} \rangle = 8 \pi^2 ({\alpha_{P,P'}})^2,$$
where
\begin{equation} \label{alpha}
\alpha_{P,P'} = \sum_i f_i(P)
\overline{f_i(P')}.
\end{equation}
\end{LEM}

\proof

Recall that on $S^2 H^{0,1}$ the scalar product is:
$$\langle a\odot b, c \odot d \rangle = \langle a, c \rangle \langle b, d
\rangle + \langle a, d \rangle \langle b, c \rangle, $$ induced by
the scalar product  $\langle a \otimes b, c \otimes d \rangle = 2
\langle a, c \rangle \langle b, d \rangle $ on $H^{0,1} \otimes
H^{0,1}$ via the inclusion of $S^2
H^{0,1}\stackrel{\iota}\hookrightarrow H^{0,1} \otimes H^{0,1}$, $
\iota(a\odot b)  = \frac{1}{2} (a \otimes b + b \otimes a).$

So, by  (\ref{csip}) one immediately computes

$$\langle \xi_P, \xi_{P'} \rangle = 8 \pi^2 \sum_{i,j} f_i(P)
\overline{f_i(P')} f_j(P) \overline{f_j(P')} =  8 \pi^2 ({\alpha_{P,P'}})^2.$$

\qed

 \section{Curvature }

 We would like now to give a formula for the curvature of the Siegel metric on $\Mg$.
 We will do the computation on the tangent vectors given by
the $\xi_P$'s. These depend on the choice of the local coordinates, but by
linearity
one can immediately derive the formulas at the tangent vectors
$\frac{\xi_P}{|\xi_P|} = \frac{\xi_P}{2 \sqrt{2} \pi \alpha_{P,P}}$, which are intrinsic.

 Recall that outside the hyperelliptic locus we have the sequence of tangent bundles:
\begin{equation}
\label{tangent} 0\rightarrow {\mathcal T}_{M_g^{(n)}} \rightarrow
j^{(n)*}{\mathcal T}_{A_g^{(n)}} \stackrel{\pi}\rightarrow
{\mathcal N} \rightarrow 0,
\end{equation}
whose dual, under the identifications  $j^{(n)*}{\Omega^1}_{A_g^{(n)}} \cong S^2 (\psi_*
K_{{\mathcal C}|{M_g^{(n)}}})$, ${\Omega^1}_{M_g^{(n)}} \cong
\psi_* ({K^2}_{{\mathcal C}|{M_g^{(n)}}})$, is

\begin{equation}
\label{dualtangent1} 0 \rightarrow {\mathcal I}_2 \rightarrow S^2
(\psi_* K_{{\mathcal C}|{M_g^{(n)}}}) \stackrel{m}\rightarrow
\psi_*({K^2}_{{\mathcal C}|{M_g^{(n)}}}) \rightarrow 0,
\end{equation}
where ${\mathcal I}_2 := {\mathcal N}^*$ and $m$ is the multiplication map.

The Hermitian connection of the variation of Hodge structures
${\mathcal R}^1 \psi_* {\Co} $, the Gauss-Manin connection,
defines a Hermitian connection on ${\mathcal F}^1= \psi_*
K_{{\mathcal C}|{M_g^{(n)}}}$, thus on ${{\mathcal F}^1}^*$, as
well as $ S^2{\mathcal F}^1$ and $S^2{{\mathcal F}^1}^*\simeq
j^{(n)*}{\mathcal T}_{A_g^{(n)}}$, which we denote by $\nabla$.

The exact sequence \eqref{tangent} defines a second
fundamental form,

$$\sigma\in Hom({\mathcal T}_{M_g^{(n)}},{\mathcal N}\otimes \Omega^1_{M_g^{(n)}}), \  \sigma: s\mapsto \pi(\nabla(s)) .$$

Similarly the exact sequence \eqref{dualtangent1} defines the
second fundamental form  $\rho \in Hom \left( {\mathcal I}_2,
\psi_*(K^2_{{\mathcal C}|{M_g^{(n)}}}) \otimes
\Omega^1_{{M_g^{(n)}}}\right)$.

The curvature form $R$ of ${\mathcal T}_{{M_g^{(n)}}}$ is computed
in terms of the curvature form $\tilde{R}$ of
${j^{(n)}}^*({\mathcal T}_{A_g^{(n)}})$ and the second fundamental
form $\sigma$. Namely, we have

\begin{equation}
\label{curvature} \langle R(s),
t\rangle =  \langle \tilde{R}(s), t \rangle -
\langle\sigma(s), \sigma(t)\rangle,
\end{equation}

where $s, t$ are local sections of  ${\mathcal
T}_{{M_g^{(n)}}}$.

At the point $[X]\in \Mg$, we need to compute
\begin{eqnarray}
\langle R(\xi_P),
\xi_{P'}\rangle(\xi_R, \overline{\xi_T}) &=& \langle
\tilde{R}(\xi_P), \xi_{P'}\rangle (\xi_R,
\overline{\xi_T}) \\ \nonumber \\ & & - \langle\sigma(\xi_P),
\sigma(\xi_{P'})\rangle(\xi_R, \overline{\xi_T}).
\end{eqnarray}

Let us now determine $\langle \tilde{R} (\xi_P), \xi_{P'}\rangle
(\xi_R, \overline{\xi_T})$ in terms of the curvature form of the
Hodge bundle.

Consider the exact sequence
\begin{equation}
\label{hodge} 0 \rightarrow {\mathcal F}^1 \rightarrow {\mathcal
R}^1 \psi_* {\Co} \otimes {\mathcal C}^{\infty}_{\Mg} \rightarrow
({\mathcal R}^1 \psi_* {\Co} \otimes {\mathcal C}^{\infty}_{\Mg})/{\mathcal
F}^1 \rightarrow 0
\end{equation}
At $[X] \in \Mg$ we have
$$0 \rightarrow H^{1,0}(X) \rightarrow H^1(X, {\Co}) \rightarrow
H^{0,1}(X) \rightarrow 0.
$$

\begin{LEM}
\label{hod}
The curvature form of the Hodge bundle is given by
$$ \langle R_{ {\mathcal F}^1}({\omega_j}),{\omega_l}\rangle(\xi_R,
\overline{\xi}_T) = 4\pi^2 \alpha_{R,T} f_j(R) \overline{f_l(T)}.
$$

\end{LEM}

\proof Since the Gauss-Manin connection on ${\mathcal R}^1 \psi_*
{\Co} \otimes {\mathcal C}^{\infty}_{\Mg} $ is flat, the following holds:
$$
\langle R_{ {\mathcal F}^1}({\omega_j}),{\omega_l}\rangle =
-\langle\epsilon(\omega_j), \epsilon(\omega_l)\rangle,
$$ where $\epsilon \in Hom(H^{1,0}(X), H^{0,1}(X) \otimes H^0(2K_X))$
is the second fundamental form of (\ref{hodge}) at the point
$[X]$. We can also view $\epsilon$ as an element in
$Hom(H^{1,0}(X) \otimes H^1(T_X), H^{0,1}(X))$, and by  a result
of Griffiths (cf. e.g. \cite{green} p.32), we have
$\epsilon(\omega_i \otimes \zeta) = \zeta(\omega_i)$. Hence, we
can write $\epsilon(\omega_j) = \sum_P (\xi_P (\omega_j) \otimes
{\xi_P}^*)$. Therefore, we have
$$
\langle\epsilon(\omega_j), \epsilon(\omega_l)\rangle = \sum_{P,S}
\langle\xi_P(\omega_j), \xi_S(\omega_l)\rangle({\xi_P}^* \otimes
\overline{{\xi_S}}^*).
$$

This implies
$$
\langle\epsilon(\omega_j), \epsilon(\omega_l)\rangle(\xi_R, \overline{\xi}_T) =
\langle\xi_R(\omega_j), \xi_T(\omega_l)\rangle.
$$

By \eqref{xdopo}, \eqref{stoka}, and \eqref{alpha} we deduce
$$ \langle \xi_R(\omega_j), \xi_T(\omega_l) \rangle= \sum_{k,i} 4\pi^2
f_j(R)f_k(R) \overline{f_l(T)f_i(T)} \langle \overline{\omega}_k ,
\overline{\omega}_i \rangle =
$$
$$
= -4 {\pi}^2 \alpha_{R,T} f_j(R) \overline{f_l(T)}.$$

\qed

Finally, we prove the following:

\begin{PROP}
\label{curvsym}
The curvature $\tilde{R}$ of  ${j^{(n)}}^*({\mathcal T}_{A_g^{(n)}}) = S^2({{\mathcal F}^1}^*)$  is given by
$$
\langle \tilde{R} (\xi_P), \xi_{P'}\rangle
(\xi_S, \overline{\xi_T}) = -64 {\pi}^4 \alpha_{S,T} \sum_{i,j,l} f_i(P) f_j(P)
\overline{f_j(P') f_l(P')} f_l(S) \overline{f_i(T)} =$$
$$= -64 {\pi}^4 \alpha_{S,T}
\alpha_{P,T} \alpha_{P, P'} \alpha_{S,P'}.$$
\end{PROP}

\proof To begin with, by Lemma \ref{csip} we have
$$\langle \tilde{R} (\xi_P), \xi_{P'}\rangle =
 4 {\pi}^2 \sum_{i,j,k,l} f_i(P) f_j(P) \overline{f_k(P')}
\overline{f_l(P')} \langle \tilde{R}(\overline{\omega_j} \odot \overline{\omega_i}),
\overline{\omega_l} \odot \overline{\omega_k}\rangle.
$$

By standard facts on complex bundles \cite{koba}, we have
\begin{equation}
\label{symm}
\langle
\tilde{R}(\overline{\omega_j} \odot
\overline{\omega_i}), \overline{\omega_l} \odot \overline{\omega_k}
\rangle
= \langle
R_{{{\mathcal F}^1}^*}(\overline{\omega_j}) \odot
\overline{\omega_i} + \overline{\omega_j} \odot R_{{{\mathcal
F}^1}^*}(\overline{\omega_i}),\overline{\omega_l} \odot
\overline{\omega_k}\rangle =
\end{equation}

$$
- \left(
{\delta}_{ik}\langle R_{{{\mathcal
F}^1}^*}(\overline{\omega_j}),\overline{\omega_l}\rangle +
{\delta}_{il}\langle R_{{{\mathcal
F}^1}^*}(\overline{\omega_j}),\overline{\omega_k}\rangle \right. + \left. {\delta}_{jl}\langle R_{{{\mathcal
F}^1}^*}(\overline{\omega_i}),\overline{\omega_k}\rangle +
{\delta}_{jk}\langle R_{{{\mathcal
F}^1}^*}(\overline{\omega_i}),\overline{\omega_l}\rangle \right). $$

Now, we observe that
$$ \langle R_{{{\mathcal
F}^1}^*}(\overline{\omega_j}),\overline{\omega_l}\rangle =\langle
R_{{\mathcal F}^1}(\omega_l), \omega_k\rangle.
$$

In fact, set $R_{{{\mathcal F}^1}^*}(\overline{\omega_j})= \sum_i
a_{ij} \overline{\omega_i}$, where $a_{ij} \in
\Omega^{1,1}_{\Mg}$. By duality, we have $R_{{{\mathcal
F}^1}}(\omega_j) = - \sum_i a_{ji} \omega_i$. Hence $ \langle
R_{{{\mathcal F}^1}^*}(\overline{\omega_j}), \overline{\omega_l}
\rangle =  -a_{lj}=\langle R_{{{\mathcal F}^1}}(\omega_l),
\omega_j \rangle. $

By \eqref{symm} and Lemma \ref{hod}, we deduce
$$
\langle \tilde{R}(\xi_P), \xi_{P'}\rangle(\xi_S,
\overline{\xi_T})= -16 {\pi}^4 {\alpha_{S,T}} \sum_{i,j,k,l}
f_i(P) f_j(P) \overline{f_k(P')} \overline{f_l(P')}\cdot$$

$$
[{\delta}_{ik}f_l(S) \overline{f_j(T)}+ {\delta}_{il} f_k(S)
\overline{f_j(T)}+ {\delta}_{jl} f_k(S) \overline{f_i(T)}+
{\delta}_{jk} f_l(S) \overline{f_i(T)}]=
$$

$$=-64 {\pi}^4 \alpha_{S,T} \sum_{i,j,l} f_i(P) f_j(P)
\overline{f_j(P') f_l(P')} f_l(S) \overline{f_i(T)}=$$
$$= -64 \pi^4 \alpha_{S,T} \alpha_{P,T} \alpha_{P,P'} \alpha_{S,P'}.$$

\qed

In order to apply \eqref{curvature}, we still need to compute
$$
\langle\sigma(\xi_P), \sigma(\xi_{P'})\rangle(\xi_S, \overline{\xi_T}).
$$

Recall that the exact sequence (\ref{tangent}) of which $\sigma$
is the second fundamental form, at $[X] \in \Mg$ is
\begin{equation}
\label{tangX} 0 \rightarrow H^1(T_X) \rightarrow S^2(H^0(K_X))^*
\rightarrow I_2(X)^* \rightarrow 0,
\end{equation}

thus
 $\sigma$
yields a homomorphism
\begin{equation}
\label{iduestella}
\sigma: H^1(T_X) \rightarrow Hom(I_2(K_X), H^0(2K_X)).
\end{equation}

Analogously, at $[X] \in \Mg$ the exact sequence
\eqref{dualtangent1} is:
\begin{equation}
\label{dualtangX} 0 \rightarrow I_2(K_X) \rightarrow S^2(H^0(K_X))
\stackrel{m}\rightarrow H^0(2K_X) \rightarrow 0,
\end{equation}
hence  the second fundamental form $\rho$ gives a homomorphism
$$ \rho: I_2(K_X) \rightarrow Hom( H^1(T_X), H^0(2K_X))$$

and for every $v \in H^1(T_X)$, and for every $Q \in I_2(X)$, we have
$$\sigma(v)(Q) = \rho(Q)(v).$$

We  recall now some results of \cite{cpt} on the second
fundamental form $\rho$. In particular we want to use Thm.2.1 and
Lemma 3.2 of \cite{cpt}, (cf. also \cite{pirola}(4.8)). Let us fix
a point $P \in X$, where $[X] \in \Mg$. We have the inclusion
$H^0(K_X(2P)) \hookrightarrow H^1(X - \{P\}, \Co) \cong H^1(X,
\Co)$. By Riemann Roch and Hodge decomposition we immediately see
that $dim (H^0(K_X(2P)) \cap H^{0,1}(X)) = 1$, so we define
$\eta_P \in H^0(K_X(2P)) \cap H^{0,1}(X)$ as the only generator of
$H^0(K_X(2P)) \cap H^{0,1}(X)$ having in a neighborhood of $P$ the
following local expression:
$$\eta_P = (-\frac{1}{(z - z(P))^2} + g(z)) dz,$$
with $g(z)$ holomorphic.

\begin{LEM} (cf. \cite{cpt} (Thm 2.1), (Lemma 3.2))
\label{CPT}
Let $Q \in I_2(K_X)$, $Q = \sum_{i,j} a_{ij} \omega_i \otimes \omega_j$, then
$$\rho(Q)(\xi_P) = -\eta_P \sum_{i,j} a_{ij} f_i(P) \omega_j \in H^0(2K_X).$$
\end{LEM}

\begin{COR}
If $X$ is any non hyperelliptic curve $ \rho$ is injective, and $
\sigma$ is non zero. In particular at any point $[X] \in \Mg$
outside the hyperelliptic locus the curvature $R$ of ${\mathcal
T}_{\Mg}$ and the curvature $\tilde{R}$ of ${j^{(n)}}^*({\mathcal
T}_{\Ag})$ are different.
\end{COR}
\proof By lemma (\ref{CPT}), for any $Q \in I_2$, $Q = \sum_{i,j}
a_{ij} \omega_i \otimes \omega_j$, $\rho(Q)(\xi_P) = 0$ implies
$\sum_{i,j} a_{ij} f_i(P) \omega_j =0,$ hence $\forall j$,
$\sum_{i} a_{ij} f_i(P) =0$. Then $Q \in ker(\rho)$ if and only if
$\sum_{i} a_{ij} f_i(P) =0$ $\forall j$, $\forall P \in X$, so
$\sum_{i} a_{ij} \omega_i =0$, which implies $Q = 0$.

Since $\sigma(\xi_P)(Q) = \rho(Q)(\xi_P)$ and $\rho$ is injective, there must exist a point $P \in X$ such that $\sigma(\xi_P) \neq 0$.
\qed
\\

Now we compute $\xi_S(\rho(Q)(\xi_P))$, where $P$ and $S$ are two
points in $X$.

Let $z$ be a local coordinate in a neighborhood of $S$, and consider a local
expression of $\rho(Q)(\xi_P) \in H^0(2K_X)$,
$$\rho(Q)(\xi_P) = \Psi_P^Q(z) dz^2.$$

\begin{LEM}
\label{Psi}
Let $Q \in I_2(K_X)$, then
$$\xi_S(\rho_Q(\xi_P)) = 2 \pi i \Psi^Q_P(S).$$

\end{LEM}

\proof
Recall that $\xi_S$ is represented by a form
$$\theta_S= \frac{1}{z} \overline{\partial} b_S \otimes
\frac{\partial}{\partial z},$$ where $z$ is a local coordinate in
a neighborhood of $S$ and $b_S$ is a bump function around $S$
which is equal to one in a neighborhood of $S$.

Let $C$ be a small circle around $S$ such that $b_S \equiv 1$ on $C$.
We have
$$\xi_S(\rho(Q)(\xi_P)) = \int_X \theta_S(\rho(Q)(\xi_P))=
\int_X \overline{\partial}\left(\frac{b_S \Psi^Q_P(z)}{z-z(S)} \right) \wedge
dz=$$
$$= \int_C \frac{\Psi^Q_P(z)}{z - z(S)} dz = 2 \pi i \Psi^Q_P(S).$$
\qed

We want now to compute $\Psi^Q_P(S).$

If $P \neq S$ the form $\eta_P$ has the following local expression in a
neighborhood of $S$:
$$\eta_P(z) = G_P(z) dz,$$
where $G_P(z)$ is holomorphic, so
$$\Psi^Q_P(z) = -G_P(z) \left( \sum_{i,j} a_{ij} f_i(P) f_j(z) \right).$$
If $P =S$, the local expression of $\eta_P$ in a neighborhood of $P$ is

$$\eta_P = (-\frac{1}{(z - z(P))^2} + g(z)) dz,$$

and we have (cf. also \cite{cpt}, Thm.3.1)
$$\Psi^Q_P(z) = -\left( -\frac{1}{(z-z(P))^2} + g(z) \right) \left( \sum_{i,j}
  a_{ij} f_i(P) f_j(z) \right)=$$
$$= \frac{ \sum_{ij} a_{ij} f_i(P) \left(f_j(P) + f'_j(P)(z - z(P)) +
    \frac{1}{2} f''_j(P) (z-z(P))^2  + h.o.t. \right)}{(z-z(P))^2} = $$
$$=  \frac{1}{2} \sum_{i,j}a_{ij} f_i(P) f''_j(P) + O(1),$$
since $Q \in I_2(K_X)$, so $\sum_{i,j} a_{ij} f_i(P)f_j(P) = 0$, and
$\sum_{i,j} a_{ij} f_i(P) f'_j(P) =0.$
Thus we have
\begin{equation}
\label{psipp} \Psi^Q_P(P) = \frac{1}{2} \sum_{i,j} a_{ij} f_i(P)
f''_j(P) = \frac{1}{2}(\mu_2(Q))(P),
\end{equation}
where $\mu_2(Q)$ is the second Gaussian map of $X$ in $Q$. For the
definition of the second Gaussian map see Section 4.

\begin{PROP}
\label{persigma}
Let $\xi_P$ be a Schiffer variation, and let $\{Q_i\}$ be an orthonormal basis
of $I_2(K_X)$, denote by $\Psi^i_P : = \Psi^{Q_i}_P$. Then the following holds:
\begin{equation}
\label{inner}
\langle\sigma(\xi_P), \sigma(\xi_{P'})\rangle(\xi_S, \overline{\xi_T}) =
4 \pi^2 \sum_i \Psi^i_P(S) \overline{\Psi^i_{P'}(T)}.
\end{equation}
\end{PROP}

\proof Fix an orthonormal basis $\{Q_i\}$ of $I_2(K_X) \subset
S^2(H^0(K_X))$. Let $\{Q_i^*\}$ be the dual basis of
$I_2(K_X)^*$. By \eqref{iduestella}, $\sigma(\xi_P) \in I_2^* \otimes
H^0(2K_X)$; hence
$$
\sigma(\xi_P) = \sum_i \sigma(\xi_P)(Q_i) \otimes Q_i^*.
$$

$\sigma(\xi_P)(Q_i) = \rho(Q_i)(\xi_P) =: \rho_{Q_i}(\xi_P) \in H^0(2K_X)$, so
$$
\sigma(\xi_P)= \sum_i
\rho_{Q_i}(\xi_P) \otimes Q_i^*.
$$

On the other hand, a basis of $H^0(2K_X)$ is given by the set $\{
\xi_S^*\}$, where $S$ runs in a set of $3g-3$ general points of
$X$. This implies that
$$
\rho_{Q_i}(\xi_P) = \sum_S
\xi_S(\rho_{Q_i}(\xi_P)) {\xi_S}^*.
$$

Therefore, the following holds:
$$
\langle\sigma(\xi_P), \sigma(\xi_{P'})\rangle(\xi_S, \overline{\xi_T}) =
$$

$$ \sum_i \sum_{V,V'} \langle\xi_V(\rho_{Q_i}(\xi_P)) {\xi_V}^*,
\xi_{V'}(\rho_{Q_i}(\xi_{P'})) {\xi_{V'}}^*\rangle(\xi_S, \overline{\xi_T}) =
$$

$$
\sum_i
\xi_S(\rho_{Q_i}(\xi_P))\overline{\xi_T(\rho_{Q_i}(\xi_{P'}))}.$$

Using lemma (\ref{Psi}) we get

$$
\langle\sigma(\xi_P), \sigma(\xi_{P'})\rangle(\xi_S, \overline{\xi_T}) =
\sum_i
\xi_S(\rho_{Q_i}(\xi_P))\overline{\xi_T(\rho_{Q_i}(\xi_{P'}))} =$$
$$=4 \pi^2 \sum_i \Psi^i_P(S) \overline{\Psi^i_{P'}(T)}.$$
\qed

From Proposition \ref{curvsym} and Proposition \ref{persigma} we
obtain a closed expression for the curvature form of ${\mathcal
T}_{\Mg}$ at $[X] \in \Mg$. More precisely, the following holds.

\begin{TEO}
\label{formula}
$$
\langle R (\xi_P), \xi_{P'}\rangle(\xi_S,
\overline{\xi_T}) = $$

$$-64 {\pi}^4 \alpha_{S,T}
\alpha_{P,T} \alpha_{P, P'} \alpha_{S,P'} - 4 \pi^2 \sum_i \Psi^i_P(S) \overline{\Psi^i_{P'}(T)}.$$

\end{TEO}

\begin{COR}
\label{holsec} The holomorphic sectional curvature of ${\mathcal
T}_{\Mg}$ at  $[X] \in \Mg$ computed at the tangent vector $\xi_P$
is given by
$$H(\xi_P) =
\frac{1}{\langle \xi_P, \xi_P \rangle \langle \xi_P, \xi_P
\rangle}\langle R(\xi_P),
\xi_P\rangle(\xi_P, \overline{\xi_P}) =
$$
$$ = -1 - \frac{1}{64 {\pi}^2 (\alpha_{P,P})^4} \sum_i |\mu_2(Q_i)(P)|^2.
$$
\end{COR}
\proof
The proof immediately follows from (\ref{formula}), (\ref{psipp}) and (\ref{scalprod}).
\qed\\

By corollary (\ref{holsec}) we see that the holomorphic sectional
curvature of $\Ag$ calculated along the tangent directions at $[X]
\in \Mg$ given by the Schiffer variations $\xi_P$ is equal to
$-1$, for all $P \in X$.

We shall now give another proof of this. We recall that the image
of the sectional curvature of $H_g$ is the segment $[-1,
-\frac{1}{g}]$ and that the tangent directions $V$ such that $H(V)
= -1$ correspond to the symmetric matrices of rank 1.

Let us now see as usual an element $\xi \in H^1(T_X)$ as a
symmetric homomorphism $H^0(K_X ) \rightarrow H^0(K_X)^*$ through
the exact sequence (\ref{tangX}). Then the above observation shows
that $H(\xi) = -1$ if and only if $\xi$ has rank one. We therefore
recall  the characterisation of the elements $\xi \in H^1(T_X)$
such that $\xi$ has rank 1.  Moreover observe that the Schiffer
variations are the points of the bicanonical curve $\phi_{2K}(X)
\subset \proj H^1(X, T_X)$. Then the statement follows as a
corollary by the following result of Griffiths and by the theorem
of Enriques-Babbage and Petri.

Define ${\mathcal X} \subset \proj H^1(X, T_X)$,
$${\mathcal X} = \{\xi \in \proj H^1(X, T_X) \ | \ rank (\xi) \leq
1\}.$$

\begin{TEO}(\cite{Gr})
Assume that $g \geq 3$ and $X$ is not hyperelliptic. Consider the image of the
bicanonical map $\phi_{2K}(X) \subset \proj H^1(X, T_X)$. Then $\phi_{2K}(X) \subset {\mathcal X}$
with equality holding if and only if the canonical curve $\phi_K(X)$ is cut
out by quadrics.
\end{TEO}

\begin{COR}
(\cite{Gr}) Assume that $g \geq 3$ and $X$ is not hyperelliptic.
Then  $\phi_{2K}(X) \subset {\mathcal X}$ with equality holding if
and only if the canonical curve $\phi_K(X)$ is not trigonal, and
it is not isomorphic to a plane quintic.

\end{COR}

\section{Second Gaussian map and holomorphic sectional curvature}

We first recall the definition of the Gaussian maps (cf.
\cite{wahl2}). Let $X$ be a smooth projective curve, $S := X
\times X$, $\Delta \subset S$ be the diagonal. Let $L$ be a line
bundle on $X$ and $L_S := p_1^*(L) \otimes p_2^*(L)$, where $p_i :
S \rightarrow X$ are the natural projections. Consider the
restriction map
$$\tilde{\mu}_{n,L}: H^0(S, L_S(-n\Delta)) \rightarrow H^0(\Delta ,
L_S(-n\Delta)_{| \Delta}).$$
Notice that since ${\mathcal O}(\Delta)_{|\Delta} \cong T_X$, we have
$$H^0(\Delta ,L_S(-n\Delta)_{| \Delta}) \cong H^0(X, 2L \otimes nK_X).$$
In the case $L = K_X$, $I_2(K_X) \subset H^0(S, K_S(-2\Delta))$,
so we can define the second Gaussian map
$$\mu_2:  I_2(K_X) \rightarrow H^0(X, 4K_X),$$
as the restriction $\tilde{\mu}_{2,K|I_2(K_X)}$.

As above we fix a basis $\{\omega_i\}$ of $H^0(K_X)$. In local coordinates
$\omega_i = f_i(z) dz$. Let $Q \in I_2(K_X)$, $Q = \sum_{i,j} a_{ij} \omega_i \otimes
\omega_j$, recall that $\sum_{i,j} a_{ij} f_if_j \equiv 0$, and since
$a_{i,j}$ are symmetric, we also have $\sum_{i,j} a_{ij} f'_if_j \equiv 0$.
The local expression of $\mu_2(Q)$ is
\begin{equation}
\mu_2(Q) = \sum_{i,j} a_{ij} f''_if_j (dz)^4= - \sum_{i,j}a_{ij} f'_i f'_j(dz)^4.
\end{equation}

We recall the following results of \cite{cfW}.

\begin{TEO}
\label{tri}(\cite{cfW} Lem.4.1, Thm.4.3) For any trigonal non
hyperelliptic curve $X$  of genus $g \geq 4$, the image of $\mu_2$
is contained in $H^0(4K_X - (q_1+...+q_{2g + 4}))$, where
$q_1+...+q_{2g+4}$ is the ramification divisor of the $g^1_3$.

 If $g \geq 8$, the rank of $\mu_2$ is $4g -18$.
\end{TEO}

We also recall
\begin{TEO}
\label{puntobase} (\cite{cfW}Thm.6.1) Assume that $X$ is smooth
curve of genus $g \geq 5$, which is non-hyperelliptic and
non-trigonal. Then for any $P \in X$ there exists a quadric $Q \in
I_2$ such that $\mu_2(Q)(P) \neq 0$. Equivalently $Im(\mu_2)\cap
H^0(4K_X-P) \neq Im(\mu_2)$, $\forall P \in X.$
\end{TEO}

Assume $[X] \in \Mg$, with $g \geq 4$, $X$ non hyperelliptic. Then corollary
(\ref{holsec}) allows us to define a function $F: X \rightarrow
\R$, given by the holomorphic sectional curvature evaluated along
the tangent vectors given by the Schiffer variations:
 $$F(P)=
H(\xi_P) = -1 - \frac{1}{64 {\pi}^2 (\alpha_{P,P})^4} \sum_i
|\mu_2(Q_i)(P)|^2 \leq -1,$$  where $\{Q_i\}$ is an orthonormal
basis of $I_2(K_X)$.

\begin{PROP}
If $g =4$, the set of points $P \in X$ such that  $F(P) = -1$ is
finite,
  which implies that $F$ is non constant.

If $g \geq 5$, $X$ not hyperelliptic, nor trigonal, then $F(P) <
-1$ for all $P \in X$.

 If $X$
is a trigonal curve of genus $\geq 4$, $F(P)= H(\xi_P)  = -1$ for every $P \in
X$ which is a ramification point of the $g^1_3$.
\end{PROP}

\proof

Assume $X$ has genus 4, then the dimension of $I_2$ is one and
$I_2$ can be generated by a quadric $Q$ of rank 4 which has norm
1. So $\forall P \in X$, $F(P) = -1 - \frac{1}{64 {\pi}^2
(\alpha_{P,P})^4}|\mu_2(Q)(P)|^2 $, hence there is a finite number
of points $P$ such that $\mu_2(Q)(P)=0$, so in these points we
have $F(P) =-1$, while $F(P) < -1$ elsewhere .

As regards the second statement, we observe that $F(P) = -1$ if
and only if $\mu_2(Q_i)(P) = 0$ for all $i$, where $\{Q_i\}$ is an
orthonormal basis of $I_2$. But then we must have $\mu_2(Q)(P) =
0$ for all $Q \in I_2$. So the proof follows by Theorem
(\ref{puntobase}).

The last statement follows from (\ref{tri}). \qed

\begin{REM}
The previous statements imply that for any curve $X \in \Mg$, not hyperelliptic, nor trigonal, for every point $P \in
X$ the holomorphic sectional curvature of $\Mg$, at
$X$ along the tangent directions given by
 $\xi_P$ is strictly smaller than the holomorphic
sectional curvature of $\Ag$. Hence the Schiffer variations are never tangent directions of totally geodesic
submanifolds of $\Ag$.

On the other hand, in the trigonal case, along
the Schiffer variations at the ramification points of the $g^1_3$, (which are a basis of the tangent space to the
trigonal locus) the
holomorphic sectional curvature of $\Mg$, coincides
with the holomorphic sectional curvature of $\Ag$.
\end{REM}

\section{The hyperelliptic locus}

We will now study the  hyperelliptic locus $HE_g \subset \Mg$. Recall that by local Torelli, the restriction of
the period map to $HE_g$ is an injective immersion (cf. \cite{os}).
Therefore we have the exact sequence
$$0 \rightarrow {\mathcal T}_{HE_g} \rightarrow {\mathcal T}_{\Ag|HE_g} \rightarrow {\mathcal N}_{HE_g| \Ag} \rightarrow 0,$$
and we denote by
$$\sigma_{HE} :  {\mathcal T}_{HE_g} \rightarrow Hom({\mathcal T}_{HE_g}, {\mathcal N}_{HE_g| \Ag})$$
the associated second fundamental form and by $\rho_{HE}$ the second fundamental form of the dual exact sequence. At the point $[X] \in HE_g$ the dual exact sequence is

$$0 \rightarrow I_2 \rightarrow S^2(H^0(K_X)) \rightarrow H^0(2K_X)^+
\rightarrow 0,$$
where $H^0(2K_X)^+$ is the invariant part of $
H^0(2K_X)$ under the hyperelliptic involution and $I_2$ is the
vector space of  the quadrics containing the rational normal
curve, so that
$$\rho_{HE}: I_2 \rightarrow Hom({\mathcal T}_{HE_g, [X]} , H^0(2K_X)^+).$$

We recall that the set of Schiffer variations at the Weierstrass
points $P_i$ generates ${\mathcal T}_{HE_g, [X]}$.

\begin{PROP}
If $X$ is hyperelliptic, $\rho_{HE}$ is injective and thus
$\sigma_{HE}$ is non zero. This implies that the curvature
$R_{HE}$ of ${\mathcal T}_{HE_g}$ is different from the curvature
$\tilde{R}$ of ${\mathcal T}_{\Ag|HE_g}$ at any point $[X] \in
HE_g$.

\end{PROP}
\proof
With the same proof of Thm 2.1, Lemma 3.2 of \cite{cpt} one can show that
$$\rho_{HE} (Q) (\xi_P) = -\eta_P \sum_{i,j} a_{ij} f_i(P) \omega_j \in H^0(2K_X)^+$$
if $P$ is a Weierstrass point of $X$ and $Q =   \sum_{i,j} a_{ij} \omega_i \otimes \omega_j \in I_2.$

So $\rho_{HE}(Q)(\xi_P) = 0$ implies $\sum_{i,j} a_{ij} f_i(P) \omega_j =0,$ hence $\forall j$, $\sum_{i,j} a_{ij} f_i(P) =0$. Then $Q \in ker(\rho_{HE})$ if and only if $\sum_{i,j} a_{ij} f_i(P) =0$ for every Weierstrass point $P \in X$. Since there are $2g+2$ Weierstrass points, this implies that $\sum_{i,j} a_{ij} \omega_i =0$, hence $Q = 0$.

Since $\sigma_{HE}(\xi_P)(Q) = \rho_{HE}(Q)(\xi_P)$ and $\rho_{HE}$ is
injective, there must exist a Weierstrass point $P \in X$ such
that $\sigma_{HE}(\xi_P) \neq 0$. \qed \\

We also observe that with the same proof as in Lemma (\ref{Psi})
and formula (\ref{psipp})  one shows that
$$\xi_P(\rho_{HE}(Q)(\xi_P))  =  \mu_2(Q)(P)$$
at a Weierstrass point $P \in X$.

Let us denote by $H_{HE}$ the holomorphic sectional curvature of
${\mathcal T}_{HE_g}$, if $[X] \in HE_g$ and  $P \in X$ is a
Weiestrass point,  we have the same expression for $H_{HE}(\xi_P)$
as in (\ref{holsec}), namely
\begin{equation}
\label{wei} H_{HE}(\xi_P) =-1 - \frac{1}{64 {\pi}^2
(\alpha_{P,P})^4} \sum_i |\mu_2(Q_i)(P)|^2
\end{equation}
where $\{Q_i\}$ is an orthonormal basis of $I_2$.

We recall now a result on the second Gaussian map proven in
\cite{cfW}.

\begin{PROP}
\label{hyp}(\cite{cfW}Lem.4.1, Prop.4.2) Let $X$ be a
hyperelliptic curve of genus $g \geq 3$. Then the rank of $\mu_2$
is $2g -5$ and its image is contained in $H^0(4K_X -
(q_1+...+q_{2g+2}))$, where $\{q_1,...,q_{2g+2}\}$ are the
Weierstrass points.
\end{PROP}

\begin{COR}
Let $[X] \in HE_g$, then $H_{HE}(\xi_P) = -1$, for any Weierstrass point $P \in X$.
\end{COR}
\proof
The proof immediately follows from (\ref{wei}) and from (\ref{hyp}).
\qed

\section{The class of the Siegel metric}

Let $\overline{M}_g$ ($\overline{\Mg})$ be the Deligne - Mumford
compactification of $M_g$ ($\Mg$). In \cite{mum1} it is shown that
the Hodge bundle extends to $\overline{M}_g$ ($\overline{\Mg})$
and its $g$-th exterior power is ample on $M_g$ ($\Mg$).

We denote by $\lambda$ both the first Chern class of the extension
of the Hodge bundle on $\overline{M_g}$ and on $\overline{\Mg}$.
We will prove that the K\"ahler form of the Siegel metric on $M_g$
extends as a closed current to $\overline{M}_g$, hence it defines
a cohomology class in $H^2(\overline{M}_g, \Co)$ which is a
multiple of $\lambda$.
\begin{TEO}
\label{class}
The K\"ahler form $\omega$ of the Siegel metric on
$M_g$ extends as a closed current to $\overline{M}_g$. Its class
$[\omega] \in H^2(\overline{M}_g, \Co)$ satisfies
$[\omega]=\pi\lambda$.
\end{TEO}

\proof On $H_g$ the Hodge metric is the only (up to multiplication
by scalars) invariant metric on the homogeneous bundle $F$.

Therefore we have an invariant metric  on the line bundle
$\Lambda^g {F}$ and thus its curvature is an invariant $(1,1)$ form $\beta$ on $H_g$.

On the other hand, the Siegel metric  is the invariant metric
obtained by the metric  on $S^2 F^*$ induced by the Hodge metric
and we denote by $\tilde{\omega}$ its K\"ahler form.

Since both $\beta$ and $\tilde{\omega}$ are invariant $(1,1)$
forms and we are on the irreducible symmetric domain $H_g$, there
exists a constant $c$ such that $\tilde{\omega} = c \beta.$ This
relation still holds on the corresponding forms on $\Ag$ which we
denote in the same way.

 In \cite{amrt} a compactification $\overline{\Ag}$ of $\Ag$ is constructed and
 it has the property that it is nonsingular and that
  $D_{\infty} := \overline{\Ag} - \Ag$ is a divisor with
  normal crossings.

In \cite{mum} it is shown that the Hodge bundle ${\mathcal F}^1$
on $\Ag$ extends as a bundle on $\overline{\Ag}$, such that the
Hodge metric has only logarithmic singularities at $D_{\infty}$.

Moreover in \cite{mum} (see also \cite{fal}), it is also proven
that the extension of the second symmetric power is isomorphic to
the sheaf of differential forms with logarithmic poles at
$D_{\infty}$:

$$S^2({\mathcal F}^1) \cong \Omega^1_{\overline{\Ag}}[D_{\infty}].$$

Furthermore Mumford proves in (\cite{mum} Thm.(3.1), Thm.(1.4))
that the extension of the Hodge metric  has ``good'' singularities
and that this implies that its first Chern class yields a closed
current on $ \overline{\Ag}$ and thus a cohomology class
$\tilde{\lambda} \in H^2(\overline{\Ag}, \Co)$.

Therefore, since on $ \Ag$ our K\"ahler form $\tilde{\omega} = c
\beta$, then also $\tilde{\omega}$ can be extended as a closed
(1,1) current on $ \overline{ \Ag}$, which we still call
$\tilde{\omega}$ and its cohomology class $[\tilde{\omega}] \in
H^2(\overline{ \Ag}, \Co)$ is given by $[\tilde{\omega}] = c
\tilde{\lambda}$.

In (\cite{nam} (18.9), see also \cite{nam1}) it is shown that the
period map $j^{(n)}: \Mg \rightarrow \Ag$ extends to a period map
$\overline{j}: \overline{ \Mg} \rightarrow \overline{\Ag}$ so we
can consider the pull-back $ \overline{j}^*([\tilde{\omega}]) \in
H^2(\overline{ \Mg}, \Co)$. Moreover, since the image of $
\overline{j}$ is not contained in the locus where the current is
singular, the pullback $\omega:=\overline{j}^*(\tilde{\omega})$ is
a well defined closed current and
$[{\omega}]=\overline{j}^*([\tilde{\omega}])$ (cf.\cite{meo}).
Moreover it gives a closed current on $\overline{M}_g$ still
denoted by $\omega$. Observe that $
\overline{j}^*(\tilde{\lambda})= \lambda$ so $[{\omega}]=c\lambda$
in $H^2(\overline{\Mg},\Co)$, hence in $H^2(\overline{M}_g,\Co)$.

In order to compute the constant $c$, we use the cycles introduced
by Wolpert in (\cite{wol}). In our case, since $[\omega]$ is a
multiple of $\lambda$, it is sufficient to compute the value of
$[\omega ]$ on the  1-dimensional family given by a varying
1-pointed elliptic curve attached to a fixed $g-1$ curve with 1
marked point ${\mathcal E}_l$ of \cite{wol}(2.2). More precisely,
let us denote by $H := \{z \in \Co \ | \ Im(z)>0\}$, by $\Gamma: =
SL(2, \Z)$, and by
$$\Gamma_l = \left\{  \left(
\begin{array}{cc}
a & b \\
c & d\\
\end{array}\right) \in \Gamma  \ | \ \left(
\begin{array}{cc}
a & b \\
c & d\\
\end{array}\right)  \equiv \left(
\begin{array}{cc}
1 & 0 \\
0 & 1\\
\end{array}\right) \ mod \ l\right\}.$$

Since the $g-1$ curve in ${\mathcal E}_l$ is constant we identify
${\mathcal E}_l$ with the curve ${H/\Gamma_l}=A_1^{(l)}$.  We have
then to compute

$$\int_{H/\Gamma_l} \omega = [\Gamma: \Gamma_l] \int_{H/\Gamma} \omega.$$
Set $E_z := \Co/(\Z \oplus z \Z)$, where $z \in H/\Gamma$, let
$\xi_z$ be a holomorphic coordinate on $E_z$, so that
$H^{1,0}(E_z) = \langle d\xi_z \rangle$. Then a cotangent
direction to the  curve ${\mathcal E}_l$ can be identified with
$d\xi_z \odot d\xi_z$ and we have: $\langle d\xi_z \odot d\xi_z,
d\xi_z \odot d\xi_z \rangle = 2 \langle d\xi_z,d\xi_z \rangle^2,$
$$\langle
d\xi_z,d\xi_z \rangle = i \int_{E_z} d\xi_z \wedge
\overline{d\xi_z} = 2 Im(z).$$ Then
$$\langle [\omega], {\mathcal E}_l\rangle  = \int_{H/\Gamma_l} \omega = i[\Gamma: \Gamma_l]
\int_D \frac{1}{8 (Im(z))^2} (dz \wedge \overline{dz}) = [\Gamma:
\Gamma_l]\frac{\pi}{12},$$ where $D$ is the fundamental domain of
the action of $\Gamma$ on $H$ and the last equality is a standard
integral computation.

Since one has $\langle \lambda, \frac{{\mathcal E}_l}{[\Gamma:
\Gamma_l]} \rangle = \frac{1}{12}$, we have
$$\frac{\pi}{12} = \langle \zeta, \frac{{\mathcal E}_l}{[\Gamma: \Gamma_l]}
\rangle =\langle c\lambda, \frac{{\mathcal E}_l}{[\Gamma:
\Gamma_l]} \rangle =  c\frac{1}{12},$$ we obtain $c = \pi$, so
finally $[\omega] = \pi \lambda.$

\qed

\vfill

\noindent {\bf Author's address:}

\bigskip

\noindent
Prof. Elisabetta Colombo\\
Dipartimento di Matematica\\
Universit\`a di Milano\\
via Saldini 50\\
     I-20133, Milano, Italy

e-mail: elisabetta.colombo@mat.unimi.it

\

\noindent
Dr. Paola Frediani\\
Dipartimento di Matematica\\
Universit\`a di Pavia\\
via Ferrata 1\\
     I-27100 Pavia, Italy

e-mail: paola.frediani@unipv.it

\end{document}